\documentstyle[12pt,leqno]{amsart}

\def\ds{\displaystyle}

\newsymbol\Subset 1362

\newcommand{\R}{{\Bbb R}}
\newcommand{\DD}{{\Bbb D}}
\newcommand{\SF}{{\Bbb S}}

\def\tl{\widetilde}

\newcommand{\loc}{{\rm loc}}

\def\es{\mathop{\rm essinf}\limits}
\newcommand{\supp}{{\rm  supp}}

\def\A{{\mathcal A}}
\def\BB{{\mathfrak B}}
\def\B{{\mathcal B}}
\def\C{{\mathcal C}}
\def\E{{\mathcal  E}}

\def\I{{\mathfrak I}}
\def\K{{\mathfrak K}}
\def\P{{\mathfrak P}}
\def\Q{{\mathcal Q}}

\def\T{{\mathcal T}}

\def\CC{{\mathfrak C}}
\def\KK{{\mathcal K}}
\def\II{{\mathcal I}}

\def\l{\ell}
\def\ba{{\bf a}}
\def\sm{\setminus}%
\def\n{{\mathbf n}}
\def\diam{{\rm diam\,}}

\newtheorem{thm}{Theorem}[subsection]  
\newtheorem{crlr}[thm]{Corollary}      
\newtheorem{defin}[thm]{Definition}    

\makeatletter
\def\theequation{\thesection.\@arabic\c@equation}
\def\thethm{\thesection.\@arabic\c@thm}
\def\thelem{\thesection.\@arabic\c@thm}
\def\thecrlr{\thesection.\@arabic\c@thm}
\def\theprp{\thesection.\@arabic\c@thm}
\def\therem{\thesection.\@arabic\c@thm}
\makeatother

\begin{document}

\title[Oblique derivative problem in generalized Morrey spaces]
{Parabolic oblique derivative problem in generalized Morrey spaces}

\author[L.G. Softova]{Lubomira G. Softova}
\address{Department of Civil Engineering,  Second University  of Naples, Via Roma 29, 81031 Aversa, Italy}
\email{luba.softova{@@}unina2.it}
\subjclass{35K20,35D35, 35B45, 35R05}
\keywords{Uniformly parabolic operator, regular oblique derivative
problem, $VMO,$  generalized Morrey spaces, unique solvability}

\maketitle

\begin{abstract}
We study the regularity of the solutions of the  oblique derivative problem for linear uniformly  parabolic equations with VMO coefficients.  We show that if the right-hand side of the parabolic equation belongs to certain generalized Morrey space $M^{p,\varphi}(Q)$ than the strong  solution belongs to the generalized 
Sobolev-Morrey  space $ \overset{\circ}{W}{}_{2,1}^{p,\varphi}(Q).$

\end{abstract}

\section{Introduction}
In the present work we consider the regular   oblique derivative problem for linear non-divergence form parabolic equations in a cylinder 
$$
\begin{cases}
 u_t- a^{ij}(x)D_{ij}u =f(x)&\quad{\rm a.e.\ in}\quad  Q,\\
 u(x',0)=0          &\quad {\rm on}\quad \Omega, \\
\partial u/\partial \l=  \l^i(x)D_iu=0	& \quad{\rm on}\quad  S.
\end{cases}
$$
The unique strong solvability of this problem was proved in \cite{Sf1}.  In \cite{Sf2}  we study the regularity of the solution in the Morrey spaces $L^{p,\lambda}$ with  $p\in(1,\infty),$ $\lambda\in(0,n+2)$ and also its H\"older regularity.  In   \cite{Sf3} we extend these studies on generalized Morrey spaces  $L^{p,\omega}$ with a weight  $\omega$ satisfying the doubling and integral conditions introduced  in \cite{Mi,  Na}. The approach associated to the names of Calder\'on and Zygmund and  developed  by Chiarenza, Frasca and Longo in \cite{CFL1, CFL2}  consists of obtaining of explicit  representation formula for the higher order  derivatives of the solution by singular and nonsingular integrals. Further  the regularity properties  of the solution follows by the continuity properties of these integrals in the corresponding spaces.  In \cite{Sf2} and \cite{Sf4} we study the regularity of the corresponding operators in the Morrey and generalized Morrey spaces while in \cite{Sf1} we dispose of the corresponding results obtained in $L^p$ by \cite{FR} and \cite{BC}.
 In recent works we study the regularity of the solutions of elliptic and parabolic problems with Dirichlet data on the boundary in Morrey-type spaces $M^{p,\varphi}$ with a weight $\varphi$ satisfying \eqref{weight}   (cf. \cite{GS1, GS2}).  Precisely,   we obtain boundedness in $M^{p,\varphi}$  for sub-linear integrals generated by singular  integrals as the Calder\'on-Zygmund. More results concerning sub-linear integrals in generalized Morrey spaces can be found in \cite{AGM, GAKS,Sf4}, see also the references therein. 

Throughout this paper the following notations are to be used, $x=(x',t)=(x'',x_n,t)\in \R^{n+1},$  $\R^{n+1}_+=\{ x'\in \R^{n},   t>0\}$ 
and $\DD_+^{n+1}=\{ x''\in \R^{n-1},   x_n>0,  t>0\},$ 
$D_iu=\partial u/ \partial x_i$, $ D_{ij}u= \partial^2 u/\partial x_i\partial
x_j$, $D_t u=u_t=\partial u/\partial t$ stand for the corresponding derivatives  while $Du=(D_1u,\ldots,D_nu)$  and  $D^2u=\{D_{ij}u\}_{i,j=1}^n$ mean the spatial gradient and the Hessian matrix of $u.$   For any measurable function $f$ and $A\subset \R^{n+1}$ we write
$$
\|f\|_{p,A}=\left(\int_A |f(y)|^pdy  \right)^{1/p},\quad f_A=\frac1{|A|}\int_A f(y)dy
$$
where $|A|$ is the Lebesgue measure of A.
Through all the paper the standard summation convention on repeated upper and lower indexes is adopted.  The letter $C$ is used for various constants and may change from one occurrence to another.

\section{Definitions and statement of the problem}\label{sec2}

Let 
$\Omega\subset \R^n,$  $n\geq 1$  be  a bounded $C^{1,1}$-domain,    
$Q=\Omega\times (0,T)$  be a cylinder in $\R^{n+1}_+,$ and  $S=\partial\Omega\times(0,T)$   
stands for the  lateral boundary of $Q.$ 
We consider the   problem
\begin{equation}\label{ODP}
\begin{cases}
{\P} u:= u_t- a^{ij}(x)D_{ij}u =f(x)&\quad{\rm a.e.\ in}\quad  Q,\\
{\I} u:= u(x',0)=0          &\quad {\rm on}\quad \Omega, \\
{\BB} u:= \partial u/\partial \l=  \l^i(x)D_iu=0	& \quad{\rm on}\quad  S.
\end{cases}
\end{equation}
under the following conditions:

\begin{itemize}
\item[$\bf (i)$] The operator  $\P$ is   supposed to be   {\it uniformly parabolic,\/}   i.e. there exists a constant 
$\Lambda>0 $ such that  for almost all (a.a.) $x\in Q$
\begin{equation}\label{2a}
\begin{cases}
&\Lambda^{-1} |\xi|^2\leq a^{ij}(x)\xi_i\xi_j\leq \Lambda |\xi|^2 \quad  \forall \xi\in\R^n,\\
&a^{ij}(x)=a^{ji}(x), \qquad i,j=1,\ldots,n.
\end{cases}
\end{equation}
The symmetry of the coefficient  matrix ${\bf a}=\{a^{ij}\}_{i,j=1}^n$
implies essential boundedness of $a^{ij}$'s and we set $\|{\bf a}\|_{\infty,Q}=
\sum_{i,j=1}^n \|a^{ij}\|_{\infty,Q}.$

\item[$\bf (ii)$] The boundary operator  $\BB$     is prescribed in terms of a
{\it directional derivative} with respect to the unit vector field
$\l(x)=(\l^1(x),\dots,\l^n(x))$   $x\in S$. We suppose that
$\BB$ is a  {\it regular oblique derivative operator\/}, i.e., the
field $    \l$ is never tangential to  $S$:
\begin{equation}\label{3a}
    \langle\l(x  )\cdot  \n  (x) \rangle =\l^i(x) \n_i(x)>0\quad  
{\rm on}\  S,\    \l^i   \in {\rm Lip}(\bar S).
\end{equation}
Here ${\rm Lip}(\bar S)$ is the class of  {\it uniformly
Lipschitz continuous functions} on $\bar S$ and
$ \n(x)$  stands for the unit outward normal to $\partial\Omega$.
\end{itemize}

In the following,  besides   the  parabolic  metric
$\varrho(x)=\max(|x'|,|t|^{1/2})$ and the defined  by it parabolic cylinders 
$$
\II_r(x)=\big\{y\in \R^{n+1}:     |x'-y'|<r, |t-\tau|<r^2\big\}\qquad |\II_r|=Cr^{n+2}\,.
$$
we  use  the equivalent  one
 $ \rho(x)=\left(\frac{|x'|^2+\sqrt{|x'|^4+4t^2}}{2}\right)^{1/2} $ (see \cite{FR}). The balls with respect to this metric are ellipsoids 
$$\E_r(x)=\big\{y\in\R^{n+1}:  \frac{|x'-y'|^2}{r^2}  +\frac{|t-\tau|^2}{r^4}<1   \big\}\qquad |\E_r|=Cr^{n+2}.$$
  Because of the equivalence of the metrics 
 all estimates obtained over ellipsoids hold true also over parabolic cylinders and in the following we shall use this   without explicit references.
\begin{defin}\label{VMO}(\cite{JN, Sa})
Let $a\in L^1_{\rm loc}(\R^{n+1}),$ denote   by 
$$
\eta_a(R)=\sup_{\E_r, r\leq R}\frac{1}{|\E_r|}\int_{\E_r}|f(y)-f_{\E_r}| dy \quad \text{ for every } R>0
$$
where $\E_r$ ranges over all ellipsoids in $\R^{n+1}.$ 
The Banach space $BMO$ (bounded mean oscillation) consists of  functions for which the following norm is finite
$$
\|a\|_{\ast}=\sup_{R>0}\eta_a(R)<\infty.
$$
A function $a$  belongs to $VMO$  (vanishing mean oscillation)  with  $VMO$-modulus $\eta_a(R)$ provided
$$
\lim_{R\to 0}\eta_a(R)=0.
 $$
For any bounded cylinder  $Q$ we define
 $BMO(Q)$ and $VMO(Q)$ taking  $a\in L^1(Q)$ and  $Q_r=Q\cap \II_r $  instead of $\E_r$ in the definition above.
\end{defin}
According to \cite{A, Jones} having a function $a\in BMO/VMO(Q)$ it is possible to extend it in the whole ${\R}^{n+1}$ preserving its
 $BMO$-norm or $VMO$-modulus, respectively. In the following we use this property without explicit references.
\begin{defin} \label{def1}
Let $\varphi:\R^{n+1}\times\R_+\to\R_+$ be a measurable function and $ p\in(1,\infty).$ A function $f\in L^p_{\loc}(\R^{n+1})$ belongs to the generalized parabolic Morrey space  $M^{p,\varphi}(\R^{n+1})$ provided
$$
\|f\|_{p,\varphi;\R^{n+1}} =\sup_{(x,r)\in\R^{n+1}\times\R_+}
\varphi(x,r)^{-1} \left(r^{-(n+2)} \, \int_{\E_r(x)}|f(y)|^p dy\right)^{1/p}<\infty.
$$
The space  $M^{p,\varphi}(Q)$  consists of  $L^p(Q)$ functions provided the following norm is finite
$$
\|f\|_{p,\varphi;Q} =\sup_{(x,r)\in Q\times\R_+}
\varphi(x,r)^{-1} \left(r^{-(n+2)} \, \int_{Q_r(x)}|f(y)|^p dy\right)^{1/p}\,.
$$
The generalized  Sobolev-Morrey  space $W_{2,1}^{p,\varphi}(Q),$ $p\in(1,\infty)$  consist of all Sobolev  functions $u\in W_{2,1}^{p}(Q)$
with distributional derivatives $D^l_tD^s_xu\in M^{p,\varphi}(Q),$  $0\leq 2l+ |s|\leq 2,$ endowed by the norm
\begin{align*}
&\|u\|_{W_{2,1}^{p,\varphi}(Q)}=\|u_t\|_{p,\varphi;Q}+ \sum_{ |s|\leq 2}\|D^s u\|_{p,\varphi;Q}.\\
& \overset{\circ}{W}{}_{2,1}^{p,\varphi}(Q)=\big\{ u\in W_{2,1}^{p,\varphi}(Q):\ \   u(x)=0 \  \   x\in
 \partial Q\big\}, \\   
&\|u\|_{\overset{\circ}{W}{}_{2,1}^{p,\varphi}(Q)}=\|u\|_{W_{2,1}^{p,\varphi}(Q)}
\end{align*}
where $\partial Q$ means the parabolic boundary $\Omega\cup\{\partial\Omega\times(0,T)\}.$
\end{defin}
\begin{thm}{\bf (Main result)}\label{main}
Let ${\bf (i)}$ and ${\bf (ii)}$ hold,  $\ba \in VMO(Q)$  and  $u\in
 \overset{\circ}{W}{}_{2,1}^{p}(Q),$ $p\in(1,\infty)$   be  a strong solution of \eqref{ODP}. If $f\in M^{p,\varphi}(Q)$ with $\varphi(x,r)$ being measurable positive function  satisfying
\begin{equation}\label{weight}
\int_{r}^{\infty} \Big(1+\ln \frac{s}{r}\Big)
\frac{\es_{s<\zeta<\infty} \varphi(x,\zeta) \zeta^{\frac{{n+2}}{p}}}{s^{\frac{{n+2}}{p}+1}}\, ds
\le C 
\end{equation}
for each $ (x,r)\in Q\times\R_+$
then $u\in \overset{\circ}{W}{}_{2,1}^{p,\varphi}(Q)$ and
\begin{equation}\label{apriori}
\|u\|_{\overset{\circ}{W}{}_{2,1}^{p,\varphi}(Q)}\leq C \|f\|_{p,\varphi;Q}
\end{equation}
with $C=C(n,p,\Lambda, \partial\Omega, T, \|\ba\|_{\infty;Q}, \eta_{\ba})$ and $\eta_{\ba}=\sum_{i,j=1}^n\eta_{a^{ij}}.$
\end{thm}
If $\varphi(x,r)= r^{(\lambda-n-2)/p}$  then $M^{p,\varphi}\equiv L^{p,\lambda}$ and the condition \eqref{weight} holds with a constant depending on $n,p$ and $\lambda.$  If $\varphi(x,r)=\omega(x,r)^{1/p}r^{-(n+2)/p} $ 
with $\omega:\R^{n+1}\times\R_+\to \R_+$ satisfying the conditions 
\begin{align*}
&k_1\leq \frac{\omega(x_0,s)}{\omega(x_0,r)}\leq k_2\qquad \forall \   x_0\in\R^{n+1}, \   r\leq s\leq 2r\\
&\int_r^\infty \frac{\omega(x_0,s)}{s}ds \leq k_3\, \omega(x_0,r)\quad k_i>0,\    i=1,2,3
\end{align*}
than we obtain the spaces $ L^{p,\omega}$  studied  in \cite{Mi, Na}. 
The following results are obtained in \cite{GS2} and treat continuity in  $M^{p,\varphi}$  of certain singular and nonsingular integrals. 
\begin{defin}\label{CZK}
A measurable function $\K(x;\xi):\R^{n+1}\times\R^{n+1}\sm\{0\}\to \R$ is called  variable parabolic Calder\'on-Zygmund kernel  (PCZK)  if:
\begin{itemize}
\item[$i)$] $\K(x;\cdot)$ is a  PCZK for a.a. $x\in\R^{n+1}:$
\begin{itemize}
\item[$a)$] $\K(x;\cdot)\in C^\infty(\R^{n+1}\sm\{0\}),$
\item[$b)$] $\K(x;\mu\xi)=\mu^{-(n+2)}\K(x;\xi)$\quad $\forall \mu>0,$
 \item[$c)$] $\ds \int_{\SF^{n}}\K(x;\xi)d\sigma_\xi=0\,,$\quad  $\ds \int_{\SF^{n}}|\K(x;\xi)|d\sigma_\xi<+\infty.$
\end{itemize}
\item[$ii)$]  $\ds \left\|D^\beta_\xi \K \right\|_{\infty;\R^{n+1}\times\SF^{n}}\leq M(\beta)<\infty$ for each  multi-index $\beta.$
\end{itemize}
\end{defin}
Consider the singular integrals 
\begin{align}\nonumber
&\KK f(x)= P.V.\int_{\R^{n+1}}\K(x;x-y)f(y)dy\\
\label{sing}
&\CC[a, f](x)= P.V.\int_{\R^{n+1}}\K(x;x-y)[a(y)-a(x)]f(y)dy\,.
\end{align}
\begin{thm}\label{CZcont}
For any $f\in M^{p,\varphi}(\R^{n+1})$ with $(p,\varphi)$ as in Theorem~\ref{main}
and $a\in BMO$   there exist  constants  depending on $n,p$ and the kernel such that
\begin{align} \nonumber
&\|\KK f\|_{p,\varphi;\R^{n+1}}
\leq C\|f\|_{p,\varphi;\R^{n+1}},\\
\label{sal22}
& \|\CC[a,f]\|_{p,\varphi;\R^{n+1}}
\leq C\|a\|_\ast\|f\|_{p,\varphi;\R^{n+1}}.
\end{align}
\end{thm}
\begin{crlr}\label{locest1}
Let $Q$ be a cylinder in $\R^{n+1}_+,$  $ f\in M^{p,\varphi}(Q),$  $a\in BMO(Q)$ and
 $\K(x,\xi):\, Q\times \R^{n+1}_+\sm \{0\}\to \R.$  Then the operators \eqref{sing} are bounded in
 $M^{p,\varphi}(Q)$ and
\begin{align}\nonumber
&\|\KK f\|_{p,\varphi;Q}\leq C\|f\|_{p,\varphi;Q},\\
\label{eq12}
&\|\CC[a,f]\|_{p,\varphi;Q}\leq C\|a\|_\ast
\|f\|_{p,\varphi;Q}
\end{align}
with $C$  independent of $a$ and $f$.
\end{crlr}
\begin{crlr} \label{locest2}
Let   $a\in VMO$ and  $(p,\varphi)$ be as in Theorem~\ref{main}.
Then for any $\varepsilon>0$ there exists a positive number
$r_0=r_0(\varepsilon,\eta_{\ba})$ such that for any  $\E_r(x_0)$
with a radius $r\in(0,r_0)$
and all $f\in M^{p,\varphi}(\E_r(x_0))$
\begin{equation}\label{normB}
\|\CC[a,f]\|_{p,\varphi;\E_r(x_0)}
\leq  C\varepsilon\|f\|_{p,\varphi;\E_r(x_0)}
\end{equation}
where $C$ is independent of $\varepsilon$, $f,$  $r$ and $x_0.$
\end{crlr}
For any $x'\in {\mathbb R}^n_+$ and any fixed $t>0$   define the {\it generalized reflection}
\begin{equation}\label{GR}
\T(x)=(\T'(x),t), \quad \T'(x) =x'-2x_n\frac{{\bf a}^n(x',t)}{a^{nn}(x',t)}
\end{equation}
where ${\bf a}^n(x)$ is  the last row of the coefficients matrix ${\bf a}(x)$ of \eqref{ODP}.
The function  $\T'(x)$ maps $\R^n_+$ into $\R^n_-$ and the kernel
$\K(x;\T(x)-y)=\K(x; \T'(x)-y',t-\tau)$ is a  nonsingular  one  for any $x,y\in \DD^{n+1}_+.$  Taking $\tl x =(x'',-x_n,t)$ there  exist positive constants $\kappa_1$ and $\kappa_2$ such that
\begin{equation}\label{CTC}
\kappa_1\rho(\tl x - y) \leq \rho({\T}(x)-y) \leq \kappa_2 \rho(\tl x -y).
\end{equation}
For any
  $f\in M^{p,\varphi}({\DD}^{n+1}_+)$ with a norm
$$
\|f\|_{p,\varphi;\DD_+^{n+1}} =\sup_{(x,r)\in\DD^{n+1}_+\times\R_+}
\varphi(x,r)^{-1} \left(r^{-(n+2)} \, \int_{\E_r(x)}|f(y)|^p dy\right)^{1/p}
$$
and  $a\in BMO({\DD}^{n+1}_+)$ define the nonsingular integral operators
\begin{align}\nonumber
&\tl   {\KK}	f(x) =\int_{{\DD}^{n+1}_+} \K (x;{\T}(x)-y)f(y) dy\\
\label{KCf}
& \tl   {\CC}  [a,f](x)=\int_{{\DD}^{n+1}_+} \K(x;{\T}(x)-y)[a(y)-a(x)]f(y) dy.
\end{align}
\begin{thm}\label{nonsing}
Let  $a \in BMO(\DD^{n+1}_+)$   and   $f\in M^{p,\varphi}(\DD^{n+1}_+) $ with   $(p,\varphi)$ as in Theorem~\ref{main}. Then the operators  $\tl\KK f$ and $\tl\CC[a, f]$ are continuous in $M^{p,\varphi}(\DD^{n+1}_+)$  and 
\begin{align}\nonumber
&\|\tl\KK f\|_{p,\varphi;\DD^{n+1}_+} \leq   C  \|f\|_{p,\varphi;\DD^{n+1}_+},\\
\label{KC}
&\|\tl\CC[a, f]\|_{p,\varphi;\DD^{n+1}_+} \leq   C \|a\|_\ast \,
\|f\|_{p,\varphi;\DD^{n+1}_+}
\end{align}
with a constant independend of $a$ and $f.$
\end{thm}
\begin{crlr} \label{localnonsing}
Let  $a\in VMO,$ then for any $\varepsilon>0$ there exists a positive number
$r_0=r_0(\varepsilon,\eta_{\ba})$ such that for any  $\E_r^+(x^0)=\E_r(x^0)\cap \DD^{n+1}_+$
with a radius $r\in(0,r_0)$ and center $x^0=(x'',0,0)$
and for  all $f\in M^{p,\varphi}(\E_r^+(x^0))$ holds
\begin{equation}\label{tlK}
\|\tl\CC[a,f]\|_{p,\varphi;\E^+_r(x^0)}
\leq  C\varepsilon\|f\|_{p,\varphi;\E^+_r(x^0)},
\end{equation}
where $C$ is independent of
$\varepsilon$, $f,$  $r$ and $x^0$.
\end{crlr}

\section{Proof of the main result}\label{sec6}

As it follows by \cite{Sf2}, the problem \eqref{ODP} is uniquely solvable in 
$ \overset{\circ}{W}{}_{2,1}^{p}(Q).$  We are going to show   that  $f\in M^{p,\varphi}(Q)$ implies 
$u\in  \overset{\circ}{W}{}_{2,1}^{p,\varphi}(Q).$
For this goal we obtain an a priori estimate of $u.$  The proof is divided in two steps.

{\it  Interior estimate.}
For any  $x_0\in \R^{n+1}_+$ consider  the parabolic semi-cylinders
$\C_r(x_0)=\B_r(x_0')\times (t_0-r^2,t_0).$  Let $v\in C_0^\infty(\C_r)$ and suppose that $v(x,t)=0$ for $t\leq 0.$ According to \cite[Theorem~1.4]{BC}  for any  $x\in \supp\, v$ the following representation formula for the second derivatives of $v$ holds true
\begin{align}
\nonumber
&D_{ij}v(x)= P.V.\int_{\R^{n+1}}\Gamma_{ij}(x;x-y)[a^{hk}(y)-a^{hk}(x)]D_{hk}v(y)dy\\
\label{RF}
&+ P.V. \int_{\R^{n+1}}\Gamma_{ij}(x;x-y)\P v(y)dy+\P v(x)\int_{\SF^n}\Gamma_j(x;y)\nu_id\sigma_y,
\end{align}
where $\nu(\nu_1,\ldots,\nu_{n+1})$ is the outward normal to $\SF^n.$  Here $\Gamma(x;\xi)$ is the fundamental solution of the operator $\P$ and $\Gamma_{ij}(x;\xi)=\partial^2
\Gamma(x;\xi)/\partial\xi_i\partial\xi_j.$  

Because of  density arguments  the representation formula  \eqref{RF}
still holds for any  $v\in W_{2,1}^{p}(\C_r(x_0)).$
 The properties of the fundamental  solution (cf. \cite{BC,LSU, Sf1}) imply  $\Gamma_{ij}$ are   Calder\'on-Zygmund kernels in the sense of Definition~\ref{CZK}.
We denote by  $\KK_{ij}$ and $ \CC_{ij}$ the singular integrals   defined in \eqref{sing} with kernels $\K(x;x-y)=\Gamma_{ij}(x;x-y).$ Then we can write  that 
\begin{align}\nonumber
D_{ij}v(x)=& \CC_{ij}[a^{hk},D_{hk}v](x)\\
\label{RF2}
& + \KK_{ij}(\P v)(x)+\P v(x)\int_{\SF^{n}} \Gamma_j(x;y)\nu_i d\sigma_y\,.
\end{align}
Because of Corollaries~\ref{locest1} and~\ref{locest2} and the equivalence of the metrics  we get
$$
\|D^2v\|_{p,\varphi;\C_r(x_0)}\leq C(\varepsilon\|D^2v\|_{p,\varphi;\C_r(x_0)}+
\|\P u\|_{p,\varphi;\C_r(x_0)}  )
$$
for some $r$  small enough. Moving  the norm of $D^2v$ on the left-hand side we get
$$
\|D^2 v\|_{p,\varphi;\C_r(x_0)}\leq C\|\P v\|_{p,\varphi;\C_r(x_0)}
$$
with a constant depending on $n,p,\eta_{\ba}(r), \|\ba\|_{\infty,Q}$ and $ \|D\Gamma\|_{\infty,Q}.$
Define  a  cut-off function $\phi(x)=\phi_1(x')\phi_2(t),$   with $\phi_1\in C_0^\infty(\B_r(x'_0)),$
$\phi_2\in C_0^\infty(\R) $ such that
$$
\phi_1(x')=\begin{cases}
1 & x'\in \B_{\theta r}(x_0')\\
0 & x'\not\in \B_{\theta'r}(x_0')
\end{cases},\qquad
\phi_2(t)=\begin{cases}
1 & t\in (t_0-(\theta r)^2, t_0]\\
0 & t< t_0-(\theta' r)^2
\end{cases}
$$
with  $\theta\in(0,1),$ $\theta'=\theta(3-\theta)/2>\theta$ and
$|D^s\phi|\leq C [\theta(1-\theta)r]^{-s},$ $s=0,1,2,$ \
$|\phi_t|\sim |D^2\phi|.$  For any solution $u\in W_{2,1}^p(Q) $ of
\eqref{ODP}
define  $v(x)=\phi(x) u(x)\in W_{2,1}^p(\C_r).$ Then we get
\begin{align*}
\|D^2u\|_{p,\varphi;\C_{\theta r}(x_0)} & \leq  \| D^2 v\|_{p,\varphi;\C_{\theta'r}(x_0)}\leq
 C\|\P v \|_{p,\varphi;\C_{\theta' r}(x_0)}\\
& \leq  C\left(\|f\|_{p,\varphi;\C_{\theta' r}(x_0)}+\frac{\|Du\|_{p,\varphi;\C_{\theta' r}(x_0)}}{\theta(1-\theta)r}
 +\frac{\|u\|_{p,\varphi;\C_{\theta' r}(x_0)}}{[\theta(1-\theta)r]^2}  \right).
\end{align*}
By  the choice   of   $\theta'$  it holds  $  \theta(1-\theta)\leq 2 \theta'(1-\theta')$ which leads to 
\begin{align*}
\big[\theta(1-\theta)r&\big]^2 \| D^2u\|_{p,\varphi;\C_{\theta r}(x_0)}\\
\leq&C \left( r^2\|f\|_{p,\varphi;Q}+\theta'(1-\theta')r\|Du \|_{p,\varphi;\C_{\theta' r}(x_0)}+
\|u\|_{p,\varphi;\C_{\theta' r}(x_0)}    \right)\,.
\end{align*}
Introducing the semi-norms
$$
\Theta_s=\sup_{0<\theta<1} \big[\theta(1-\theta)r \big]^s \|D^s u \|_{p,\varphi;\C_{\theta r}(x_0)}\qquad s=0,1,2
$$
and taking the supremo with respect to $\theta$ and  $\theta'$  we get
\begin{equation}\label{theta}
 \Theta_2
 \leq C
\left(r^2\|f\|_{p,\varphi;Q} +\Theta_1+\Theta_0 \right)\,.
\end{equation}
The  interpolation inequality  \cite[Lemma~4.2]{Sf3} gives that there  exists a positive constant $C$  independent of $r$ such that
$$
\Theta_1\leq \varepsilon\, \Theta_2+\frac{C}{\varepsilon}\, \Theta_0\qquad \text{ for any } \varepsilon\in(0,2).
$$
Thus   \eqref{theta} becomes
$$
[\theta(1-\theta)r]^2\|D^2u\|_{p,\varphi;\C_{\theta r}(x_0)}\leq\Theta_2\leq
 C\left(r^2\|f\|_{p,\varphi;Q}+\Theta_0 \right)
$$
for each  $ \theta\in(0,1).$
Taking  $\theta =1/2$ we get  the    Caccioppoli-type estimate
$$
\|D^2u\|_{p,\varphi; \C_{r/2}(x_0)} \leq  C\left( \|f\|_{p,\varphi;Q}+
 \frac1{r^2}\|u\|_{p,\varphi;\C_r(x_0)} \right).
 $$
To estimate  $u_t $ we exploit the parabolic structure of  the equation
and the boundedness of the coefficients
\begin{align*}
\|u_t\|_{p,\varphi;\C_{r/2}(x_0)}& \leq \|{\bf a}\|_{\infty;Q} \|D^2u\|_{p,\varphi;\C_{r/2}(x_0)}
+\|f\|_{p,\varphi;\C_{r/2}(x_0)}\\
&  \leq C\big( \|f\|_{p,\varphi;Q} + \frac1{r^2}\|u\|_{p,\varphi;\C_r(x_0)} \big).
\end{align*}
 Consider  cylinders $Q'=\Omega'\times(0,T)$ and $Q''=\Omega''\times(0,T)$ with $\Omega'\Subset\Omega''\Subset\Omega,$
by standard covering procedure and partition of the unity we get
\begin{equation}\label{intest}
\|u\|_{W_{2,1}^{p,\varphi}(Q')}  \leq C
     \big(\|f\|_{p,\varphi;Q}
+\|u\|_{p,\varphi;Q''}\big)
\end{equation}
where  $C$ depends on $n, p, \Lambda, T,
\|D\Gamma\|_{\infty;Q}, \eta_{\ba}(r),$ $\|{\bf a}\|_{\infty,Q}$ and $\text{dist}(\Omega',\partial\Omega'').$

{\it Boundary estimates.}
For any fixed $R>0$ and $x^0=(x'',0,0)$  define the  semi-cylinders
$$
\C_R^+(x^0) =\C_R(x^0)\cap  \DD^{n+1}_+\,.
 $$
Without lost of generality we can take $x^0=(0,0,0).$ Define\\
 $\B_R^+=\{|x'|<R, x_n>0\},$
$S_R^+=\{|x''|<R, x_n=0, t\in(0,R^2)\}$ and consider the problem
\begin{equation}\label{VP}
\begin{cases}
{ \P} u  := u_t- a^{ij}(x)D_{ij}u =f(x) & {\rm a.e.\	in}\quad
\C_R^+,\\
 {\I} u := u(x',0)=0 & {\rm on} \quad  \B_R^+,\\
 {\BB} u:= \l^i(x)D_i u=0     & {\rm on} \quad S_R^+.
\end{cases}
\end{equation}
Let $u\in W^{p}_{2,1}(C_R^+)$ with  $u=0$ for $t\leq 0$ and $x_n\leq 0,$ then the following {\it  representation formula\/}   holds  (see  \cite{MPS,Sf1})
$$D_{ij}u(x)=I_{ij}(x)-J_{ij}(x)+H_{ij}(x)$$
 where
\begin{align*}
I_{ij}(x)=& P.V.
\int_{\C_R^+ }\Gamma_{ij}(x;x-y)F(x;y) dy\\
&+f(x)
\int_{\SF^{n}}\Gamma_j(x;y)\nu_i
d\sigma_y,\qquad  i,j=1,\dots,n\,;
\end{align*}
\begin{align*}
J_{ij}(x)=& \int_{\C_R^+}\Gamma_{ij}(x;\T(x)-y)F(x;y) dy;\\
J_{in}(x)=& J_{ni}(x)=\int_{\C_R^+}\Gamma_{il}(x;\T(x)-y)
\left(\frac{\partial \T(x)}{\partial x_n} \right)^l F(x;y)dy,\\
&\qquad\qquad\qquad\qquad\qquad\qquad\qquad i,j=1,\ldots,n-1\\
J_{nn}(x)=& \int_{\C_R^+} \Gamma_{ls}(x;\T(x)-y)
\left( \frac{\partial \T(x)}{\partial x_n} \right)^l
\left(\frac{\partial \T(x)}{\partial x_n} \right)^s
F(x;y)dy;
\end{align*}
\begin{align*}
F(x;y)=&  f(y)+[a^{hk}(y)-a^{hk}(x)]D_{hk}u(y)\\
 H_{ij}(x)=& (G_{ij}\ast_2 g)(x) +g(x'',t)\int_{\SF^n} G_j(x;y'',x_n,\tau) n_i d\sigma_{(y'',\tau)},\\
&\qquad\qquad\qquad\qquad\qquad\qquad\qquad   i,j=1,\ldots,n,\\
\frac{\partial T(x)}{\partial x_n}=& \left(-2\frac{a^{n1}(x)}{a^{nn}(x)},\ldots,-2\frac{a^{n n-1}(x)}{a^{nn}(x)},-1  \right).
\end{align*}
Here the kernel $G=\Gamma \Q,$ is a byproduct of the fundamental solution and  a  bounded regular  function $\Q.$ Hence its derivatives $G_{ij}$ behave as $\Gamma_{ij}$ and the convolution that appears in $H_{ij}$ is defined as follows 
\begin{align*}
&(G_{ij}\ast_2 g)(x) = P.V.
\int_{S_R^+}G_{ij}(x;x''-y'',x_n,t-\tau)g(y'',0,\tau) dy''d\tau,\\
&g(x'',0,t)= \left[\big(\l^k(0)-\l^k(x'',0,t)\big)D_ku
- \l^k(0) (\Gamma_k\ast F)\right]\Big\vert_{x_n=0}(x'',0,t),\\
&(\Gamma_k\ast F)(x)= \int_{\C_R^+} \Gamma_k(x;x-y)F(x;y) dy\,.
\end{align*}
Here $I_{ij}$ are a sum of  singular  integrals and bounded surface integrals hence the estimates obtained in Corollaries~\ref{locest1} and~\ref{locest2} hold true.
 On the nonsingular integrals $J_{ij}$ we apply the estimates obtained in Theorem~\ref{nonsing} and Corollary~\ref{localnonsing} that give
\begin{equation}\label{IJ}
\|I_{ij}\|_{p,\varphi;\C_R^+}+  \|J_{ij}\|_{p,\varphi;\C_R^+}\leq C\big(\|f\|_{p,\varphi;\C_R^+}+\eta_{\ba}(R)\|D^2u\|_{p,\varphi;\C_R^+}   \big)
\end{equation}
for all  $i,j=1,\ldots,n.$ 
To estimate the   norm of $H_{ij}$ 
we suppose that the vector
field $\l$  is extended in $\C_R^+$  preserving its
Lipschitz regularity and the norm. This automatically leads	 to extension  of
the function  $g$ in $ \C_R^+$ that is
\begin{equation}\label{gg}
g(x)=\big(\l^k(0)-\l^k(x)\big)D_ku(x)- \l^k(0)	(\Gamma_k\ast F)(x)\,.
\end{equation}
Applying the estimates for the heat potentials   \cite[Chapter~4]{LSU} and the trace theorems in $L^p$ \cite[Theorems~7.48,~7.53]{Ad} (see also  \cite[Theorem~1]{Sf1}) we get
$$
\int_{ \C_R^+}  |(G_{ij}\ast_2 g)(y)|^p dy \leq C
\left(\int_{\C_R^+}|g(y)|^pdy + \int_{\C_R^+}|Dg(y)|^p dy\right)\,.
$$
Taking   a parabolic cylinder $\II_r(x)$ centered in some point $x\in \C_R^+$ we have
\begin{align*}
\int_{ \C_R^+\cap\II_r(x)}  |(G_{ij}\ast_2 g)(y)|^p dy& \leq C  \frac{r^{n+2}}{\varphi(x,r)^{-p}} 
\Big(\frac{\varphi(x,r)^{-p}}{ r^{n+2}}\int_{\C_R^+\cap \II_r(x)}|g(y)|^pdy\\
& +\frac{\varphi(x,r)^{-p}}{ r^{n+2}}
\int_{\C_R^+\cap \II_r(x)}|Dg(y)|^p dy\Big)\\
& \leq C  \frac{\varphi(x,r)^{-p}}{ r^{n+2}}
\left(\|g\|^p_{p,\varphi;\C_R^+} +\|Dg\|^p_{p,\varphi;\C_R^+}\right)\,.
\end{align*}
Moving $ \frac{\varphi(x,r)^{-p}}{ r^{n+2}}$  on the left-hand side and 
  taking the supremo with respect to $(x,r)\in \C_R^+\times\R_+$	we get
$$
 \|G_{ij}\ast_2 g\|^p_{p,\varphi;\C_R^+}\leq C\left(\|g\|^p_{p,\varphi;\C_R^+}
+\|Dg\|^p_{p,\varphi;\C_R^+} \right).
$$
An immediate consequence of \eqref{gg}	 is the estimate
\begin{align*}
\|g\|_{p,\varphi;\C_R^+}& \leq \|[\l^k(0)-\l^k(\cdot) ]D_ku\|_{p,\varphi;\C_R^+}    +\|  \l^k(0)	(\Gamma_k\ast F) \|_{p,\varphi;\C_R^+}\\
&    \leq
CR\|\l\|_{{\rm Lip}(\bar S)}\|D u\|_{p,\varphi;\C_R^+} + \|\Gamma_k\ast f\|_{p,\varphi;\C_R^+}\\
&+\| \Gamma_k\ast[a^{hk}(\cdot)- a^{hk}(x)]D_{hk}u \|_{p,\varphi;\C_R^+}.
\end{align*}
The convolution $\Gamma_k\ast f$  is  a   Riesz
potential. On the other hand 
\begin{align*}
|(\Gamma_k\ast f)(x)|\leq& C  \int_{\C_R^+}\frac{|f(y)|}{\rho(x-y)^{n+1}}\,dy\\
&\leq
 CR \int_{\C_R^+}\frac{|f(y)|}{\rho(x-y)^{n+2}}\,dy\leq  C  \int_{\C_R^+}\frac{|f(y)|}{\rho(x-y)^{n+2}}\,dy
 \end{align*}
with a constant depending on $T$ and $\diam\Omega.$ Because of   \cite[Lemma~7.12]{GT} 
$$
\| \Gamma_k\ast f\|_{p,\C_R^+}\leq C\|f\|_{p,\C_R^+}
$$
which allows  to apply \cite[Theorem~3.3]{GS2} that gives
$$
\|\Gamma_k\ast f\|_{p,\varphi;\C_R^+} \leq C \|f\|_{p,\varphi;\C_R^+}.
$$
Analogously
$$
|\Gamma_k\ast [a^{hk}(\cdot)-a^{hk}(x)]D_{hk}u(\cdot)|\leq
C\int_{\C_R^+}\frac{|a^{hk}(y)-a^{hk}(x)||D_{hk}u(y)|}{\rho(x-y)^{n+2}}dy
$$
with a constant depending on $\diam\Omega$ and $T.$
The kernel  $\rho(x-y)^{-(n+2)}$ is a nonnegative singular one and  the \cite[Theorem~0.1]{B} gives
$$
\|\Gamma_k\ast  [a^{hk}(\cdot)-a^{hk}]D_{hk}u\|_{p,\C_R^+} \leq C \|\ba\|_\ast \|D^2u\|_{p;\C_R^+}.
$$
Applying again the results for sub-linear integrals \cite[Theorem~3.7]{GS2} we get
$$
\|\Gamma_k\ast  [a^{hk}(\cdot)-a^{hk}]D_{hk}u\|_{p,\varphi;\C_R^+} \leq C \|\ba\|_\ast 
\|D^2u\|_{p,\varphi;\C_R^+}\,.
$$
Hence 
\begin{align}\label{g}
\|g\|_{p,\varphi;\C_R^+}\leq&
C \big( R \|\l\|_{{\rm Lip}(\bar S)}\|D u\|_{p,\varphi;\C_R^+}+\|f\|_{p,\varphi;\C_R^+}\\
\nonumber
&+   
R\eta_{\ba}(R)\|D^2 u\|_{p,\varphi;\C_R^+} \big).
\end{align}
Further, the Rademacher theorem asserts existence   almost everywhere 
of the derivatives $D_h \l^k\in L^\infty$, thus
$$
D_h g(x)=-D_h\l^k(x) D_ku(x) + [\l^k(0) - \l^k(x)]D_{kh}u -\l^k(0)(\Gamma_{kh}\ast F)(x).
$$
The $M^{p,\varphi}$ norm of the last term  is estimated as above  and
\begin{align}\label{Dg}
\nonumber
\|D g\|_{p,\varphi;\C_R^+}\leq& C\big(\|D\l\|_{\infty;S}\|D u\|_{p,\varphi;\C_R^+}+ 
R\|\l\|_{{\rm Lip}(\bar S)}\|D^2u\|_{p,\varphi;\C_R^+}\\
&+ \|f\|_{p,\varphi;\C_R^+} + \eta_{\ba}(R)\|D^2 u\|_{p,\varphi;\C_R^+}\big).
\end{align}
Finally unifying \eqref{IJ}, \eqref{g} and \eqref{Dg} we get
\begin{align*}
\|D^2u\|_{p,\varphi;\C_R^+} \leq &  C\Big( \|f\|_{p,\varphi;Q} +(1+R) \|Du\|_{p,\varphi;\C_R^+} \\
&+ (R+\eta_{\ba}(R)+R\eta_{\ba}(R))\|D^2u\|_{p,\varphi;\C_R^+} \Big)
\end{align*}
with a constant depending on known quantities and $\|\l\|_{{\rm Lip}(\bar S)}$ and
 $ \|D\l\|_{\infty;S}. $
Direct calculations lead to an interpolation inequality in $M^{p,\varphi}$ analogous to \cite[Lemma~3.3]{LSU} (cf. \cite{Sf3})   
$$
\|Du\|_{p,\varphi;\C_R^+}\leq \delta\|D^2u\|_{p,\varphi;\C_R^+}+\frac{C}{\delta} 
\|u\|_{p,\varphi;\C_R^+},\quad \delta\in(0,R)\,.
$$
Taking $0<\delta=\frac{R}{R+1}<R$ we get 
\begin{align*}
\|D^2u\|_{p,\varphi;\C_R^+} \leq &  C\Big( \|f\|_{p,\varphi;Q} +R \|D^2u\|_{p,\varphi;\C_R^+} +
\frac{C}{R}\|u\|_{p,\varphi;\C_R^+}\\
&+ (R+\eta_{\ba}(R)+R\eta_{\ba}(R))\|D^2u\|_{p,\varphi;\C_R^+} \Big)\,.
\end{align*}
Choosing $R$ small enough end moving the terms containing the norm of $D^2u$ on the
 left-hand side we get
$$
\|D^2u\|_{p,\varphi;\C_R^+}\leq C\left(\|f\|_{p,\varphi;\C_R^+}+\frac1{R}\|u\|_{p,\varphi;\C_R^+}   \right)\,.
$$
Because of the parabolic structure of the  equation analogous estimate holds also for $u_t.$
Further the {\it Jensen inequality}  applied to   $u(x)=\int_0^t u_s(x',s)ds$    gives 
 $$\|u\|_{p,\varphi;\C_R^+} \leq C R^2\|u_t\|_{p,\varphi;\C_R^+}\leq	C\left(R^2
\|f\|_{p,\varphi;\C_R^+}+R\|u\|_{p,\varphi;\C_R^+}\right).
$$
Choosing $R$  smaller, if necessary,  we get $\|u\|_{p,\varphi;\C_R^+}
\leq
C\|f\|_{p,\varphi;\C_R^+}$ and therefore
\begin{equation}\label{LB}
\|u\|_{W^{2,1}_{p,\varphi}(\C_R^+)} \leq C\|f\|_{p,\varphi;\C_R^+}\leq
C\|f\|_{p,\varphi;\C_R^+}.
 \end{equation}
Making a covering $\{\C_\alpha^+\},$ $\alpha\in{\A}$ such that
  $Q\setminus Q'\subset \bigcup_{\alpha\in{\A}} \C_\alpha^+,$
considering
a partition of  unity  subordinated to that  covering and applying 
\eqref{LB}   for each $\C^+_\alpha$  we get
\begin{equation}\label{boundest}
\|u\|_{W^{2,1}_{p,\varphi}(Q\setminus Q')} \leq C\|f\|_{p,\varphi;Q} 
\end{equation}
with a constant depending on $ n,$ $p,$ $\Lambda,$ $ T,$  $ \diam\Omega,$
 $\|D\Gamma\|_{\infty;Q},$  $ \eta_{\ba},$  $\|{\bf a}\|_{\infty;Q},$  $ \|\l\|_{{\rm Lip}(\bar S)},$ and $  \|D\l\|_{\infty,S}.$

The estimate \eqref{apriori}  follows from \eqref{intest}  and \eqref{boundest}.

\end{document}